
%
%
\font\sc=cmcsc10
\font\cmrXVII=cmr17
\font\cmrVIII=cmr8
\font\cmbxIX=cmbx9
\font\slVIII=cmsl8

%
\def\cstar{$C^*$}
\def\ep{\varepsilon}
\def\epp{{\ep^\prime}}
\def\eppn{{\ep^\prime_n}}
\def\square{\hbox{$\sqcap\!\!\!\!\sqcup$}}
\def\e{\hbox{\rm e}}
\def\ker{\hbox{\rm Ker}}
\def\dist{\hbox{\rm dist}}
\def\det{\hbox{\rm det}}
\def\Z{{\bf Z}}
\def\phibar{\overline{\phi}}
\def\for{\hbox{\quad for \quad}}
\def\un{u(n)}
\def\sun{su(n)}
\def\S1{{\bf S}^1}
\def\up{u^\prime}
\def\vp{v^\prime}
  \def\MRClass{\cmrVIII 1991 \slVIII MR Subject Classification:
  \cmrVIII
  46L05, 
  49R20.}
  \def\support{\cmrVIII Partially supported by FAPESP, Brazil. On leave
from the University of S\~ao Paulo.}

%
\def\section#1{\bigskip \bigskip \noindent {\bf #1} \medskip}
\def\statement#1{\bigskip \noindent {\bf #1} \sl}
\def\proof{\medskip \noindent {\bf Proof.} \rm}
\def\endproof{\hfill \square}
\def\ref#1#2#3#4{\medskip \item {[#1]}{\sc #2}, {\sl #3}, #4.}

%
\def\choi    {{\bf 1}}
\def\dixmier {{\bf 2}}
\def\soft    {{\bf 3}}
\def\fell    {{\bf 4}}
\def\rieffel {{\bf 5}}
\def\voic    {{\bf 6}}

%
\magnification=1200
\baselineskip=14pt

%
\centerline {\cmrXVII THE SOFT TORUS II}

\vskip 2truemm
\centerline {\cmbxIX A Variational Analysis of}

\centerline {\cmbxIX    Commutator Norms}

\vskip 7truemm
\centerline{Ruy Exel\footnote{$\null^\dagger$}{\support}
\footnote{\null}{\MRClass}}

\vskip 7truemm
\baselineskip=9pt
\parshape 1 15truemm 135truemm
\noindent
{\cmbxIX Abstract. \cmrVIII The field of \cstar-algebras over the
interval [0,2] for which the fibers are the Soft Tori is shown to be
continuous.  This result is applied to show that any pair of
non-commuting unitary operators can be perturbed (in a weak sense) in
such a way to decrease the commutator norm. Perturbations in norm are
also considered and a characterization is given for pairs of unitary
operators which are local minimum points for the commutator norm in
the finite dimensional case.}

\baselineskip=12pt
\section{1. Introduction.}

As in [\soft], for every $\ep$ in the real interval [0,2] we let
$A_\ep$ be the universal unital \cstar-algebra generated by unitary
elements $u_\ep$ and $v_\ep$ subject to the relation $$|| u_\ep v_\ep
- v_\ep u_\ep|| \leq \ep.$$

Clearly, if $\ep_1 \leq \ep_2$ there is a unique homomorphism
$A_{\ep_2 } \longrightarrow A_{\ep_1}$ sending $u_{\ep_2 }$ and
$v_{\ep_2 }$ respectively to $u_{\ep_1}$ and $ v_{\ep_1}$.  In case
$\ep_2 = 2$ and $\ep_1 = \ep$ we shall denote this map by $\phi_\ep$,
i.e.  $$\phi_\ep : A_2 \longrightarrow A_\ep$$

One of the main results of the present work (Theorem 3.4) is the fact
that there exists a continuous field of \cstar-algebras over the
interval [0,2] such that $A_\ep$ is the fiber over $\ep$ and moreover
such that $$\ep \in [0,2] \mapsto \phi_\ep (a) \in A_\ep$$ is a
continuous section for every $a$ in $A_2$. We refer the reader to
[\dixmier] for a treatment of the theory of continuous fields of
$C^*$-algebras.

The central point in proving our main result is to show that
$||\phi_\ep (a)||$ is a continuous function of $\ep$ for all $a$ in
$A_2$.  If we let $$J_\ep = \ker(\phi_\ep)$$ then we have that
$||\phi_\ep (a) || = \dist(a, J_\ep)$ and, since the $J_\ep$'s clearly
form a \underbar{decreasing} chain of ideals, $|| \phi_\ep (a) ||$ is
seen to be an \underbar{increasing} function of $\ep$.

Let us denote by $J_\ep^+$ the closure of the union of all $J_\epp$
for $\epp > \ep$ and by $J_\ep^-$ the intersection of all $J_\epp$ for
${\epp < \ep}$. That is $$J_\ep^+ = \overline{\bigcup_{\epp > \ep}
J_\epp} $$ and $$J_\ep^- = \bigcap_{\epp < \ep} J_\epp.$$

\statement{1.1 Proposition.}
Let $\ep$ be in [0,2) (resp. (0,2]). A sufficient condition for
$$\epp \longrightarrow || \phi_{\epp} (a) ||$$
to be right (resp. left) continuous at $\ep$ for all $a$ in
$A_2$ is that $J_\ep = J_\ep^+$ (resp. $J_\ep = J_\ep^-$). \rm

\proof Note that for $a$ in $A_2$ we have $$\dist( a, J_\ep^+ ) =
\inf_{\epp > \ep} \dist (a, J_{\epp})$$ and that $$\dist( a, J_\ep^- )
= \sup_{{\epp < \ep}} \dist (a, J_{\epp}).$$

The first identity follows from trivial Banach space facts. On the
other hand, the second one cannot be generalized to Banach spaces so
let's prove it.

Let $$\phi: A_2 \longrightarrow \prod_{{\epp < \ep}} A_\epp$$ be given
by $\phi(a) = (\phi_\epp (a))_{{\epp < \ep}}$ and observe that $\ker(
\phi ) = J_\ep^-$ so $$ \dist (a, J_\ep^- ) = || \phi (a) || =
\sup_{{\epp < \ep}} || \phi_\epp (a) || = \sup_{{\epp < \ep}} \dist
(a, J_\epp).$$ The conclusion now follows without much trouble.
\endproof

\bigskip  En passant, let's give an example to support our statement that
the above fact does not generalize to Banach spaces.

Let $E = C[0,1]$ and for $f$ in $E$ put $$||f|| = |f(0)| + \sup_{t \in
[0,1]} |f(t)|.$$ Let $$E_n = \{ f \in E: f([1/n ,1]) = 0\}$$ and let
$g$ be the constant function $g = 1$.  Then $$\dist(g, E_n) = 1$$ for
all $n$ while $$\dist(g, \bigcap_{n \geq 1} E_n) = || g ||= 2.$$

\statement{1.2 Proposition.} For all $\ep \in [0,2)$ one has $J^+_\ep
= J_\ep$. \rm

\proof If we denote by $u^+_\ep$ and
$v^+_\ep$ the images of $u_2$ and $v_2$ in $A_2 /J^+_\ep$
then it is clear that
$$|| u^+_\ep v^+_\ep - v^+_\ep u^+_\ep || \leq\ep$$
hence by the universal property of $A_\ep$ there is a homomorphism
$$A_\ep = A_2 /J_\ep \longrightarrow A_2 / J_\ep^+ $$
sending $u_\ep$ and $v_\ep$ to $u^+_\ep$ and $v^+_\ep$.  Therefore
$J_\ep \subseteq J_\ep^+$. But since the reverse inclusion is trivial,
our proof is complete. \endproof

\bigskip The main technical result of this work, which shall be proven
in the next two sections, is the following.

\statement{1.3 Theorem.}
For all $\ep$ in $(0,2]$ one has $J_\ep^- = J_\ep$. \rm

\bigskip Therefore we have

\statement{1.4 Corollary.}
For all $a$ in $A_2$ the function
$$\ep \in [0, 2] \longrightarrow || \phi_\ep (a) ||$$
is continuous. \rm

\section{2. The Case $\ep < 2$.}

Throughout this section we shall fix a real number $\ep$ with
$0 < \ep < 2.$

\statement{2.1 Lemma.}
Suppose $u_0, u_1, ... , u_n$ is a finite sequence of unitary
elements in a \cstar-algebra $B$ such that $|| u_{k-1} - u_k || \leq
\ep$ for all $k$.  Then, for every $\delta > 0$, there are
unitaries $v_0, v_1, ... , v_n$ in $B$ such that

\medskip
(i) $|| v_k - u_k || \leq \delta$ \for $k = 0, ... , n$

\medskip
(ii) $||v_{k-1} - v_k || < \ep$ \for $k= 1, ... , n.$ \rm

\proof Since $|| u_{k-1} - u_k || \leq  \ep
< 2$, it follows that $|| u_k u_{k-1}^{-1} - 1 || < 2$
so that $-1$ is not in the spectrum of $u_k u_{k-1}^{-1}$.
Therefore we may define
$$h_k = \log (u_k u^{-1}_{k-1} )$$
where $\log$ is the principal branch of the logarithm.  Each
$h_k$ is then a skew adjoint element in $B$ and we have
$$u_k = \e^{h_k} u_{k-1} \hbox{\quad for \quad} k = 1, 2, ... , n.$$

We shall choose our sequence $v_0, ... , v_n$ of the form
$$v_0 = u_0$$
and
$$v_k = \e^{t_k h_k} u_{k-1} \hbox{\quad for \quad} k \geq 1$$
where each $t_k$ will be a suitably chosen positive real
number approaching 1 from below.

     A real function which will be useful in our estimates is
$$d(x) = | 1 - \e^{ix} | = 2 \sin \left( {x \over 2} \right)
\for x \in \hbox{\bf R}.$$

For instance, observe that if $h$ is skew adjoint and $|| h ||
\leq \pi$ one has
$$|| 1 - \e^h || = d ( || h ||)$$
by the spectral theorem. Moreover, for all $k = 1, ... , n$
$$d (|| h_k ||) = || 1 - \e^{h_k} || = || u_{k-1} - u_k || \leq \ep$$
which implies (since $d$ is increasing in [0, $\pi$]) that
$$|| h_k || \leq \theta$$
where $\theta = d^{-1} (\ep)$. Note that $\theta < \pi$ because $\ep <
2$.

In search of the correct choice of the $t_k$'s observe that
$$|| v_0 - v_1 || = || u_0 - \e^{t_1  h_1 } u_0 || =
d (t_1 || h_1 ||) \leq d (t_1 \theta)$$
and that for $k \geq 1$ we have
$$|| v_k - v_{k+1} || \leq || v_k - u_k || + || u_k - v_{k+1} || = $$
$$||\e^{t_k h_k} u_{k-1} - \e^{h_k} u_{k-1} || +
|| u_k - \e^{t_{k+1}h_{k+1}} u_k || = $$
$$d ((1-t_k)|| h_k||) + d (t_{k+1} || h_{k+1} ||) \leq $$
$$d ( (1-t_k) \theta) + d (t_{k+1} \theta).$$

Note that $d^\prime (x) = \cos (x/2)$ so for $x \in [0, \theta]$
we have $m \leq d^\prime (x) \leq 1$
where $m = \cos (\theta / 2)$ is strictly positive since
$\theta < \pi$.  By the mean value theorem we then have
$$m | t - s | \leq | d (t) - d (s) | \leq | t - s |$$
for all $t$ and $s$ in $[0, \theta]$.

If this last fact is used in our previous computations, we obtain
$$|| v_0 - v_1 || \leq d (t_1 \theta) =
\ep - (d (\theta) - d (t_1 \theta)) \leq $$
$$\ep - m (\theta - t_1 \theta) = \ep - m\theta (1 - t_1)$$
while for $k \geq 1$
$$|| v_k - v_{k+1} || \leq d ( ( 1 - t_k) \theta) +
d (t_{k+1} \theta) \leq$$
$$( 1 -t_k) \theta + \ep - (d (\theta) - d (t_{k+1}\theta)) \leq$$
$$(1 - t_k) \theta + \ep - m(\theta - t_{k+1} \theta) = $$
$$\ep + (1 -t_k -m (1- t_{k+1})) \theta.$$

Therefore the condition that $|| v_{k-1} - v_k || < \ep$
will be fulfilled as long as
$$1 - t_1 > 0$$
and
$$1 - t_{k+1} > {{1 - t_k} \over m}.$$

If we thus put $t_k = 1 - 2^k \sigma / m^k$
for $\sigma < (m/2)^n$ we have that each $t_k$ is in (0,1) and
$|| v_{k-1} - v_k || < \ep$.

Clearly, as $\sigma$ tends to zero, each $v_k$ approaches the
corresponding $u_k$ so a suitable choice for $\sigma$ yields
$$|| v_k - u_k || \leq \delta$$
for all $k = 0, 1, ... , n.$ \endproof

\bigskip Recall from [\soft] that $A_\ep$ is isomorphic to the crossed
product
$$A_\ep \simeq B_\ep \times_\tau \Z$$
where $B_\ep$ is the universal \cstar-algebra generated
by a sequence $\{ u^{(\ep)}_n : n \in \Z \}$ of
unitaries satisfying the relations
$$|| u^{(\ep)}_n - u^{(\ep)}_{n+1} || \leq \ep$$
for all $n$.

Moreover $\tau$ is the automorphism of $B_\ep$ given by
$$\tau (u^{(\ep)}_n ) = u^{(\ep)}_{n+1} \for n \in \Z.$$

\statement{2.2 Proposition.}
There exists a sequence $(\psi_n)_{n \in {\bf N}}$ of
endomorphisms of $B_\ep$, converging pointwise to the identity map,
such that
$$ \sup_{k \in \Z} || \psi_n (u^{(\ep)}_k) -
\psi_n (u^{(\ep)}_{k+1} ) || < \ep.$$ \rm

\proof By the previous Lemma, let for every $n$
$$v^{(n)}_{-n}, v^{(n)}_{-n+1}, ... , v^{(n)}_{0}, v^{(n)}_1, ... ,
v^{(n)}_n$$
be unitaries in $B_\ep$ such that
$$|| v^{(n)}_k - u^{(\ep)}_k || \leq {1 \over n} \for k = -n , ... ,n$$
and
$$|| v^{(n)}_{k-1} - v^{(n)}_k || < \ep \for k = -n+1, ... , n.$$

      Define $\psi_n: B_\ep \longrightarrow B_\ep $ by
$$\psi_n ( u^{(\ep)}_k ) = \left\{
\matrix{
v^{(n)}_{-n}   & \hbox{if} & k< -n \cr
\null \cr
v^{(n)}_k      & \hbox{if} & -n \leq k \leq n \cr
\null \cr
v^{(n)}_n      & \hbox{if} & k > n.} \right.$$

It is then clear that
$$\lim_{n \rightarrow \infty} \psi_n (u^{(\ep)}_k ) =
u^{(\ep)}_k$$
which implies that $\psi_n$ converges pointwise to the identity.

The condition that
$$\sup_{k \in \Z} || \psi_n (u^{(\ep)}_k) - \psi_n (u^{(\ep)}_{k+1})
|| < \ep$$
is also clearly satisfied. \endproof

\statement{2.3 Definition.}
Let $K_\ep$ be the ideal in $B_2$ given by the kernel of the
canonical map
$$\phi_\ep: B_2 \longrightarrow B_\ep.$$

\statement{2.4 Theorem.}
One has \quad $\bigcap_{{\epp < \ep}} K_\epp = K_\ep$. \rm

\proof Let $x$ be in \quad $\bigcap_{\epp < \ep} K_\epp$ and put
$y = \psi_\ep (x)$.  We have
$$y = \lim_{n\rightarrow \infty} \psi_n (y) = \lim_{n\rightarrow
\infty} \psi_n (\phi_\ep (x)).$$

Now let
$$\eppn = \sup_k || \psi_n \phi_\ep (u^{(2)}_k)
- \psi_n \phi_\ep (u^{(2)}_{k+1}) || = \sup_k || \psi_n
(u^{(\ep)}_k ) - \psi_n (u^{(\ep)}_{k+1})||$$
which by (2.2) is
strictly less than $\ep$.  So $\psi_n \phi_\ep$ factors through
$B_\eppn$ and since $x$ is in $K_\eppn$ we have $\psi_n
\phi_\ep (x) = 0$.  So $y = 0$ which implies that $x$ is in
$K_\ep$. The converse inclusion is trivial. \endproof

\statement{2.5 Lemma.}
Let $\phi: A \longrightarrow B$ be a \cstar-algebra homomorphism
and suppose $\phi$ is equivariant with respect to automorphisms $\alpha$
and $\beta$ of $A$ and $B$ respectively.

Let $K = \ker (\phi)$ and $J = \ker (\phibar)$
where $\phibar$ is the canonical extension of $\phi$
to the corresponding crossed products by \Z.  Then
$$J = \{  x \in A \times_\alpha \Z: E_A (xu^{-n}) \in K
\for n \in \Z \} $$
where
$$E_A  :A \times_\alpha \Z \longrightarrow A$$
is the associated conditional expectation [\rieffel] and $u$ is the
unitary implementing $\alpha$. \rm

\proof It is clear that \  $\phi \circ E_A = E_B \circ
\phibar$ \ where $E_B$ is the conditional expectation for $B
\times_\beta \Z$.

     Let $v$ be the implementing unitary for $B \times_\beta \Z$.
Given $x \in A \times_\alpha \Z$ we have that $x$ is in $J$ if and
only if $\phibar(x) = 0$ which is equivalent to the fact that
$E_B(\phibar(x)v^{-n}) = 0$ for all $n$, or that $\phi(E_A(xu^{-n}) =
0$. But this is to say that $E_A(xu^{-n})$ is in $K$. \endproof

\statement{2.6 Theorem.}
For every $\ep \in (0,2)$ one has $J_\ep^- = J_\ep.$ \rm

\proof Let $E: A_2 \longrightarrow B_2$ be the
conditional expectation induced by the isomorphism
$$A_2 \simeq B_2 \times_\tau \Z.$$

Given
$x$ in $J_\ep^-$ we have that $E(x u^{-n})$ is in $K_\epp$
for all $n$ and $\epp < \ep$ so
$$E(xu^{-n}) \in \bigcap_{\epp < \ep} K_\epp = K_\ep$$
for all $n$, which shows that $x \in J_\ep$. The converse inclusion is
trivial. \endproof

\section{3. The case $\ep = 2$.}

The purpose of this section is to prove that $J_2^- = J_2$ or,
since $J_2 = (0)$, that
$$\bigcap_{\ep < 2} J_\ep = (0).$$

Note that the techniques employed in the previous section do not
work here because one couldn't take logarithms in the proof of
(2.1) if $\ep = 2$.  A different approach is thus necessary.

\statement{3.1 Lemma.}
Let $w_1$ and $w_2$ be $n \times n$ unitary matrices.  Then
$|| w_1 - w_2 || = 2$ if and only if $\det (w_1  + w_2) = 0$. \rm

\proof We have that $|| w_1 - w_2 || =
2$ if and only if $|| w_1 w_2^{-1} -1 || = 2$ which is
equivalent to $- 1$ being in the spectrum of $w_1
w_2^{-1}$ which is to say that $\det (w_1 w_2^{-1} + 1)
= 0$ or that $\det (w_1  + w_2) = 0$.  \endproof

\statement{3.2 Proposition.}
Given unitary $n \times n$ matrices $u$ and $v$ such that
$|| uv - vu || = 2$
there is, for every $\delta > 0$, a unitary $u^\prime$ with
$$|| u^\prime - u || < \delta$$
and
$$|| u^\prime v - v u^\prime|| < 2.$$ \rm

\proof Write $u = \e^h$ for some skew adjoint $h$
and let $u(t) = u\e^{-th}$ for all real $t$.  Put
$$f(t) = \det \left( u(t) v + v u(t) \right)$$
and observe that
$$f(1) = \det (2v) \not = 0$$
while
$$f(0) = \det (u v + v u) = 0$$
by (3.1).  Therefore $f$ is not a constant function and
since it is analytic, its zeros are isolated.  So there are
arbitrarily small values of $t$ for which $f(t) \not = 0$, which
is to say, by (3.1) again, that
$$|| u (t) v - v u (t) || < 2 .$$

Taking $t$ sufficiently small will also ensure that
$|| u (t) - u || < \delta$. \endproof

\statement{3.3 Theorem.}
One has $J^-_2 = J_2$. \rm

\proof Assume by way of contradiction that $a \in J^-_2$ is
non-zero.

Note that $A_2$ is isomorphic to the full \cstar-algebra of the
free group on two generators so that by [\choi] $A_2$ has a separating
family of finite dimensional representations.  Let therefore
$$\pi : A_2 \longrightarrow M_n ({\bf C})$$
be a representation such that $\pi (a) \not = 0$.

Let $u = \pi (u_2) $ and $v = \pi (v_2) $ and write $$u =
\lim_{i \rightarrow \infty} u^\prime_i$$
where $|| u^\prime_i v - v u^\prime_i || < 2$ by (3.2).  For each $i$
let $\pi_i $ be the representation of $A_2$ such that
$$\pi_i (u_2) = u^\prime_i \hbox{\quad and \quad} \pi_i (v_2)  = v.$$

Then, if $\ep_i = || u^\prime_i v - v u^\prime_i ||$, we
have that $\pi_i$ vanishes on $J_{\ep_{i}}$ and thus also on
$J^-_2$ so $\pi_i (a) = 0$.

On the other hand it is clear that $\pi_i$ converges pointwise
to $\pi$ so $\pi (a) = 0$ which is a contradiction. \endproof

\statement{3.4 Theorem.} There exists a continuous field of
\cstar-algebras over the interval [0,2] such that $A_\ep$ is the fiber
over $\ep$ and such that $$\ep \in [0,2] \mapsto \phi_\ep (a) \in
A_\ep$$ is a continuous section for every $a$ in $A_2$.

\proof Let $S$ be the set of sections $$\ep \in [0,2] \mapsto \phi_\ep
(a) \in A_\ep$$ for $a$ in $A_2$.  According to [\dixmier]
(Propositions 10.2.3 and 10.3.2) all one needs to check is that $S$ is
a *-subalgebra of the algebra of all sections, that the set of all
$s(\ep)$ as $s$ runs through $S$ is dense in $A_\ep$ and that
$||s(\ep)||$ is continuous as a function of $\ep$ for all $s$ in $S$.

The first two properties are trivial while the last one follows from
(1.1), (1.2), (2.6) and (3.3). \endproof

\section{4. Local Minima for Commutator Norms.}

As clearly indicated by the results obtained above, the phenomena
under consideration is related to the following question

\statement{4.1 Question.}
Given unitary operators $u$ and $v$ which do not commute, when
is it possible to perturb $u$ and $v$ in order to obtain a new pair
$u^\prime$ and $v^\prime$ such that
$$|| u^\prime v^\prime - v^\prime u^\prime || < || uv - vu ||? $$ \rm

\bigskip In other words (4.1) calls for a characterization of pairs of
unitary operators which are not local minimum points for the commutator
norm.  Proposition (3.2) is clearly a partial answer to (4.1) and it
says that when $||uv-vu||=2$ in finite dimensions then $(u,v)$ is
never such a local minimum point.

In formulating the question above we chose not to specify the precise
meaning of ``to perturb'' in order to allow for different points of view.

The following is a complete answer to (4.1) under quite a loose type
of perturbation which we might call {\it *-strong dilated
perturbation}.

\statement{4.2 Theorem.}
Let $u$ and $v$ be non commuting unitary operators on a Hilbert
space $H$.  Then there are nets $(u_i)_i$ and $(v_i)_i$ of
unitary operators on $H^\infty$ (the direct sum of infinitely many
copies of $H$) satisfying
$$|| u_i v_i - v_i u_i || < || uv - vu ||$$
and such that the compressions
$$\hbox{proj}_H \circ u_i|_H \hbox{\quad and \quad}
\hbox{proj}_H \circ v_i |_H$$
of $u_i$ and $v_i$ to $H$ converge *-strongly to $u$ and $v$
respectively. \rm

\proof After Theorems (2.6) and (3.3) this basically becomes a
consequence of [\fell] and some well known results on representation
theory.  We therefore restrict ourselves to a sketch of the proof.

Let $\ep = || u v - v  u ||$ and consider the set $N$ of
states on $A_2$ which vanish on some $J_\epp$ for $\epp < \ep$.
Since $\bigcap_{\epp < \ep} J_{\epp} = J_\ep$ it can be proved that $N$ is
weakly dense in the set of states of $A_2$ that vanish on $J_\ep$.

Given $u$ and $v$ let $\pi$ be the representation of $A_2$ on $H$
such that $\pi (u_2) = u$ and $\pi (v_2) = v$.  Assume without
loss of generality that $\pi$ is cyclic with cyclic vector $\xi$ and put
$$f: a \in A_2  \longrightarrow \langle \pi (a) \xi, \xi \rangle .$$

Since  $\pi $ factors through $A_\ep$ we have that $f$
vanish on $J_\ep$ so there exists a net $(f_i)_{i\in I}$ in $N$
converging weakly to $f$.  Let, for every $i, \pi_i$ be the GNS
representation of $A_2$ corresponding to $f_i$.
Since each $f_i$ vanish on some $J_{\ep^\prime_i}$ with
$\ep^\prime_i < \ep$ the same is true for $\pi_i$ hence
$$|| \pi_i (u_2) \pi_i (v_2) - \pi_i (v_2) \pi_i (u_2) || \leq
\ep^\prime_i < \ep.$$

     Using the methods of [\fell] one may assume that the
space $H_i$ where $\pi_i$ acts is a subspace of $H^\infty$ and that
the conclusion holds with
$$u_i = \pi_i (u_2) + 1 - p_i$$
and
$$v_i = \pi_i (u_2) + 1 - p_i$$
where $p_i$ is the projection onto $H_i$. \endproof

\bigskip We therefore see that there are no local minimum points for the
commutator norm other than the commuting pairs, as long as we consider
*-strong dilated perturbations.

The situation is quite different if norm perturbations are considered
as we shall see in the next Section.

\section{5. Norm Perturbations in Finite Dimensions.}

Let us now study question (4.1) for pairs of unitary operators on a
finite dimensional Hilbert space. From now on we shall only consider
norm perturbations.

For $n \geq 3$ denote by $\Omega_n$ and $S_n$ the $n\times n$
Voiculescu's  unitary matrices (see [\voic])
$$\Omega_n = \pmatrix{
\omega &\        &\        &\      &\        \cr
\      &\omega^2 &\        &\      &\        \cr
\      &\        &\omega^3 &\      &\        \cr
\      &\        &\        &\ddots &\        \cr
\      &\        &\        &\      &\omega^n \cr
}
\hbox{\quad and \quad}
S_n = \pmatrix {
0  &\  &\      &\  &1  \cr
1  &0  &\      &\  &\  \cr
\  &1  &0      &\  &\  \cr
\  &\  &\ddots &\  &\  \cr
\  &\  &\      &1  &0  \cr
}$$
where $\omega = \e^{2 \pi i/n}$.

\statement{5.1 Theorem.}
For $n \geq 3$ there exists a neighborhood $V$ of the pair
$(\Omega_n, S_n)$ in $U (n) \times U(n)$ so that
$$|| uv- vu || \geq || \Omega_n S_n - S_n \Omega_n ||$$
for all pairs $(u, v)$ in $V$.

\proof Note that $\Omega_n S_n \Omega^{-1}_n
S^{-1}_n = \omega I_n$ so that if $(u,v)$ is close enough to
$(\Omega_n, S_n)$ then the spectrum of $uv u^{-1} v^{-1}$ is in a
small neighborhood of $\omega$ in the complex plane.

On the other hand, notice that
$$\det (u v u^{-1} v^{-1}) = 1$$
so that if the spectrum of $uv u^{-1} v^{-1}$ is the set
$\{ \e^{i\theta_1 } , ... , \e^{i\theta_n}\}$
with $-\pi < \theta_i < \pi$ one has that
$\sum^n_{k=1} \theta_k$ is in $2\pi \Z.$
So, by continuity,
$$\sum^n_{k=1} \theta_k = 2\pi.$$

Therefore, for some $k_0$ we must have $\theta_{k_0} \geq
2 \pi / n$ and it follows that
$$|| uv - vu || \geq | \e^{i\theta_{k_0}} - 1 |  \geq | \omega -1 | =
|| \Omega_n S_n - S_n \Omega_n || .$$ \endproof

\bigskip In other words, the pair $(\Omega_n, S_n)$ is a local minimum
for the commutator norm. This shows that the situation is quite different
from what we saw when we considered *-strong dilated perturbations
in Section (4). This result should also be compared with [\voic].

Clearly the method used in Theorem (5.1) above applies to show that
the conclusion is also true for any pair of unitary matrices whose
multiplicative commutator is a scalar multiple of the identity, but
not equal to $-I$ (see 3.2).

In fact, among irreducible pairs there are no other examples as we
shall prove shortly. We say that a pair of unitary operators is
irreducible when there is no proper invariant subspace for both
elements of the pair.

Denote by $\gamma$ the map
$$\gamma: (u, v) \in U (n) \times U (n) \mapsto u v u^{-1}
v^{-1} \in SU (n).$$

\statement{5.2 Lemma.}
A point $(u,v) \in U (n) \times U (n)$ is regular for $\gamma$ (in the
sense that $\gamma$ is a submersion at $(u,v)$) if and only if $(u,v)$
is an irreducible pair. \rm

\proof If $h,k$ are in the Lie algebra $\un$ of $U(n)$, a simple
computation shows that
$$d\gamma_{(u,v)} (uh, vk) = uv (v^{-1} hv - h + k - u^{-1} ku)
u^{-1} v^{-1} .$$

So $\gamma$ is a submersion at $(u,v)$ if and only if the map
$$L: (h, k) \in \un \times \un \longrightarrow
v^{-1} h v -h + k - u^{-1} ku \in \sun$$
is onto the Lie algebra $\sun$ of $SU(n)$.

Under the inner product on $\sun$ defined by
$$\langle x, y \rangle\ =\ \hbox{Trace}(xy^*)$$
the orthogonal space to the image of $L$ can easily be seen to be the set
$$\{ x \in \sun: xu = ux \hbox{\ and \ } xv = vx \}.$$

Now, by Schur's lemma, irreducibility of $(u,v)$ can be
characterized by the fact that only scalars commute with both $u$ and
$v$.  Since $\sun$ contains no scalar matrices the
result follows.  \endproof

\statement{5.3 Theorem.}
If $(u,v)$ is an irreducible pair in $U(n) \times U(n)$ and at
the same time a local minimum for the commutator norm then $u v
u^{-1} v^{-1}$ is a scalar. \rm

\proof By the open mapping theorem the image
under $\gamma$ of a neighborhood of $(u, v)$ is a neighborhood of
$\gamma (u, v)$.  But since
$$||uv - vu||  = || \gamma (u, v) - 1 ||$$
it follows that $\gamma (u,v)$ is a local minimum for the map
$$w \in SU(n) \longrightarrow || w - 1 ||$$
and this implies, as a moments thought will reveal, that
$\gamma (u,v)$ is a scalar.  \endproof

\bigskip This completes the classification of local minima for
irreducible pairs.  So let us now consider a reducible pair $(u,v)$ of
unitary $n \times n$ matrices.  As usual write
$$u = \oplus u_j   \quad \hbox{and}\quad v = \oplus v_j$$
where each pair $(u_j, v_j)$ is irreducible and observe that
$$|| uv - vu || = \max_j || u_j v_j - v_j u_j ||.$$

\statement{5.4 Theorem.}
Let $u = \oplus u_j$ and $v = \oplus v_j$ be as above and suppose that
the pair $(u,v)$ is a local minimum for the commutator norm.  Then
for some value of $j$ for which
$$|| u_j v_j - v_j u_j || = || uv - uv ||$$
one has that $u_j v_j u_j^{-1} v^{-1}_j$ is a scalar. \rm

\proof If this is not so then for all such $j$
the pair $(u_j, v_j)$ admits by (5.3) a small perturbation decreasing
the commutator norm.  Together, these perturbations yield a
contradiction to the hypothesis. \endproof

\bigskip A natural question which one could ask is, of course, whether the
converse to Theorem (5.4) is also true.  A good test case is given by
the pair $$(\Omega_n \oplus I_m, S_n \oplus I_m)$$ that is, the direct
sum of Voiculescu's unitaries with the $m \times m$ identity matrix.

This pair clearly satisfies the conclusion of (5.4) and so it is
natural to ask whether or not it is a local minimum for the commutator
norm.

Despite strong favorable evidence given by some partial positive
results and a large amount of computer simulation supporting this
thesis, we were unable to establish a proof for this fact. In fact we
do not even know whether the above pair is a local minimum for
$n=3$ and $m=1$. Nevertheless, we conjecture that

\statement{5.5 Conjecture.}
The converse of (5.4) is also true. \rm

\section{6. An Example.}

Considering the apparent discrepancy between (4.2) and (5.1) it is
perhaps interesting to see a concrete example of nets $(u_i)_i$ and
$(v_i)_i$, whose existence is guaranteed by (4.2), in case the given
unitaries are taken to be Voiculescu's unitaries, i.e. $u=\Omega_n$
and $v=S_n$.

For that purpose it is enough to find, for all $\delta > 0$, unitary
operators $\up$ and $\vp$ on a separable, infinite dimensional Hilbert
space $H$, and an orthonormal set $\{\xi_k: k \in \Z/n\Z\}$ of vectors
in $H$ such that

\medskip
(i) $||\up \vp - \vp \up|| < ||\Omega_n S_n - S_n \Omega_n||$

(ii) $||\up (\xi_k) - \omega^k \xi_k|| < \delta$ \quad and

(iii) $||\vp(\xi_k) - \xi_{k+1}|| < \delta$

\medskip
\noindent for all $k$ in $\Z/n\Z$ where $\omega = \e^{2 \pi i/n}$.

Let $H=L_2(\S1)$ and let $\up$ be the unitary operator on $H$ defined by
$$\up(\xi)|_z = z\xi(z) \for \xi \in H, \ z \in \S1$$ and, for $\theta
< 2\pi/n$, let $\vp$ be defined by $$\vp(\xi)|_z =
\xi(\e^{-i\theta}z) \for \xi \in H, \ z \in \S1.$$

Let $\xi_0$ be a unit vector in $H$ represented by a function $f$ on
$\S1$ supported in a neighborhood $V$ of $z=1$ which is small enough
so that $V$ is disjoint from $\e^{ik\theta} V$ for $k=1, 2, ..., n-1$.
For all such $k$ let $\xi_k = {\vp}^k(\xi_0)$. The reader may now
check that (i), (ii) and (iii) hold as long as $\theta$ is close to
$2\pi/n$ and the diameter of $V$ is small enough.

\bigskip
\centerline{\sc References}

\medskip
\ref{\choi}
{M. D. Choi}
{The full $C^*$-algebra of the free group on two generators}
{Pacific J. Math.  {\bf 87}(1980), 41-48}

\ref{\dixmier}
{J. Dixmier}
{$C^*$-Algebras}
{North Holland, 1982}

\ref{\soft}
{R. Exel}
{The Soft Torus and applications to almost commuting matrices}
{Pacific J. Math, to appear}

\ref{\fell}
{J. M. G. Fell}
{$C^*$-algebras with smooth dual}
{Illinois J. Math.  {\bf 4}(1960), 221-230}

\ref{\rieffel}
{M. A. Rieffel}
{Induced representations of \cstar-algebras}
{Advances in Math. {\bf 13}(1974), 176-257}

\ref{\voic}
{D. Voiculescu}
{Asymptotically commuting finite rank unitary operators without
commuting approximants}
{Acta Sci. Math. (Szeged) {\bf 45}(1983), 429-431}

\vskip 20truemm
\parindent=0truemm
\sc
Current address:

Department of Mathematics and Statistics

University of New Mexico

Albuquerque, NM 87131, USA

\rm
e-mail:  exel@math.unm.edu

\vskip 15truemm
\sc
Permanent address:

Departamento de Matem\'atica,

Universidade de S\~ao Paulo,

Caixa Postal 20570,

01498 S\~ao Paulo SP, Brasil
\end